\documentclass[12pt]{amsart}

\usepackage{graphics}
\usepackage{amsmath}
\usepackage{amsfonts}
\usepackage{amssymb}
\usepackage{amscd}

\newtheorem{theorem}{\bf Theorem}

\title{Siciak-Zahariuta extremal functions \\ and polynomial hulls}

\author{Finnur L\'arusson}
\address{School of Mathematical Sciences, University of Adelaide, Adelaide SA 5005, Australia.} 
\email{finnur.larusson@adelaide.edu.au}

\author{Ragnar Sigurdsson}
\address{Science Institute, University of Iceland, Dunhaga 3, IS-107 Reykjav\'ik, Iceland.} 
\email{ragnar@hi.is}


\thanks{The first-named author was supported in part by the Natural Sciences and Engineering Research Council of Canada.}

\dedicatory{Dedicated to Professor J\'ozef Siciak on the occasion of his 75th birthday}

\date{8 January 2007}

\begin{document}

\begin{abstract} 
We use our disc formula for the Siciak-Zahariuta extremal function to characterize the polynomial hull of a connected compact subset of complex affine space in terms of analytic discs.  
\end{abstract}

\maketitle

\noindent
The polynomial hull $\hat K$ of a compact subset $K$ of complex affine space $\mathbb C^n$ is the compact set of those $x\in\mathbb C^n$ for which $|P(x)|\leq \sup_K |P|$ for all complex polynomials $P$ in $n$ variables.  The set $K$ is said to be polynomially convex if $\hat K=K$.  Polynomial hulls are usually difficult to determine.  By the maximum principle for holomorphic functions it is clear that $x\in\hat K$ if there is a continuous map $f$ from the closed unit disc $\overline{\mathbb D}$ into $\mathbb C^n$, holomorphic on $\mathbb D$, such that $f(0)=x$ and $f$ maps the unit circle into $K$.  In the early days of polynomial convexity theory it seemed possible that $\hat K$ might simply be the union of all analytic discs in $\mathbb C^n$ with boundary in $K$.  Stolzenberg's counterexample of 1963 \cite{Stolzenberg} showed that $\hat K\setminus K$ can be nonempty---and in fact quite large: see \cite{DuvalLevenberg}---without containing any nonconstant analytic discs.  The question of whether polynomial hulls could nevertheless be somehow described in terms of analytic discs remained open for three decades, until Poletsky derived an answer from his theory of disc functionals \cite{Poletsky} (see Theorem \ref{poletsky} below).  In this note, we give a different characterization of the polynomial hull of a connected compact subset of $\mathbb C^n$, based on our generalization in \cite{LarussonSigurdsson2} of Lempert's disc formula for the Siciak-Zahariuta extremal function from the convex case to the connected case.  The gist of both Poletsky's result and ours is to suitably weaken the requirement that the analytic discs in the description of the polynomial hull map the whole unit circle continuously into the compact set.

We start by reviewing some preliminaries.  Recall that a holomorphic function $f$ on the unit disc $\mathbb D$ belongs to the Nevanlinna class $N$ if $\log|f|$ has a positive harmonic majorant on $\mathbb D$.  Equivalently, $f$ is a quotient $g/h$, where $g$ and $h$ are bounded holomorphic functions on $\mathbb D$ and $h$ has no zeros.  Then $f$ has nontangential boundary values $f^*(\zeta)$ at almost every point $\zeta$ of the boundary $\mathbb T$ of $\mathbb D$, and the measurable function $\log|f^*|$ is integrable with respect to the normalized arc length measure $\sigma$ on $\mathbb T$.  Unless $f$ is identically zero, $f$ factors uniquely into a product of a Blaschke product, an outer function, and a singular function $s$, given by the formula
$$s(z)=\exp\int_{\mathbb T}\frac{\zeta+z}{\zeta-z}\,d\mu(\zeta),$$
where $\mu$ is a finite real measure on $\mathbb T$, singular with respect to $\sigma$.  If $\mu$ is positive, then $s$ is the reciprocal of an inner function.  The least harmonic majorant of $\log|f|$ is the Poisson integral of the finite real measure $\log|f^*|d\sigma+d\mu$.  We call $\mu$ the singular measure of $f$. 

We say that an analytic disc $f:\mathbb D\to\mathbb C^n$ is Nevanlinna if each of its components is a Nevanlinna function or, equivalently, the subharmonic function $\log\|f\|$ has a positive harmonic majorant on $\mathbb D$.  Here, $\|\cdot\|$ denotes the Euclidean norm (although any other norm on $\mathbb C^n$ would do). 

For a thorough account of Nevanlinna functions, see \cite{Garnett}, Sec.\ II.5.  

\begin{theorem}  Let $K\subset\mathbb C^n$ be compact and connected.  For $a\in\mathbb C^n$, the following are equivalent.
\begin{enumerate}
\item[(i)]  $a$ is in the polynomial hull $\hat K$ of $K$.
\item[(ii)]  For every neighbourhood $U$ of $K$, there is a Nevanlinna disc $f:\mathbb D\to\mathbb C^n$ with $f(0)=a$ and nontangential boundary values in $U$ almost everywhere on $\mathbb T$, whose components have positive, arbitrarily small singular measures.
\end{enumerate}
\label{characterization}
\end{theorem}

\begin{proof}  (ii) $\Rightarrow$ (i):  Let $P$ be a nonconstant complex polynomial in $n$ variables.  Let $U$ be a neighbourhood of $K$ and let $f$ be a Nevanlinna disc as in (ii).  It is easily seen that the Nevanlinna class is a $\mathbb C$-algebra, so $P\circ f$ is Nevanlinna.  If we denote by $\mu$ the singular measure of $P\circ f$ and by $h$ the least harmonic majorant of $\log|P\circ f|$, then
$$\log|P(a)| \leq h(0) = \int_{\mathbb T}\log|P\circ f^*|d\sigma +\mu(\mathbb T) \leq \sup_U \log|P| + \mu(\mathbb T).$$
It remains to show that $\mu(\mathbb T)$ can be made smaller than any preassigned positive number.  This holds, since for reasons explained in the next paragraph, there is an integer $m$, only depending on $P$, such that $\frac 1 m \mu\leq \mu_1+\dots+\mu_n$, where $\mu_j$ is the singular measure of the $j$-th component of $f$.  

Clearly, when Nevanlinna functions are multiplied, their singular measures get added.  For $j=1,2$, let $g_j$ be a Nevanlinna function with singular measure $\mu_j$.  Let $h_j$ be the smallest positive harmonic majorant of $\log|g_j|$, that is, the smallest harmonic majorant of $\log^+|g_j|$.  Then $h_j$ is the Poisson integral of $\log^+|g_j^*|d\sigma+d\mu_j^+$.  Crude estimation shows that $|g_1+g_2|\leq e^{h_1}+e^{h_2}\leq 2e^{h_1+h_2}$, so $h_1 + h_2 +\log 2$ is a positive harmonic majorant for $\log|g_1+g_2|$.  It follows that the singular measure of $g_1+g_2$ is at most $\mu_1^+ + \mu_2^+$.

Note that so far, we have only used the compactness of $K$.

(i) $\Rightarrow$ (ii):  It is well known that $a\in\hat K$ if and only if $V_K(a)=0$, where $V_K$ is the (unregularized) Siciak-Zahariuta extremal function of $K$.  Furthermore, $V_K=\sup V_U$, where $U$ runs through any basis of neighbourhoods of $K$ in $\mathbb C^n$ (see \cite{Klimek}, 5.1.5).  Since $K$ is connected, these neighbourhoods may be taken to be connected.  It follows that $a\in\hat K$ if and only if $V_U(a)=0$ for every connected open neighbourhood $U$ of $K$.

Now assume $a\in\hat K$ and let $U$ be a neighbourhood of $K$.  By shrinking $U$ if necessary, we may assume that $U$ is open and connected.  Let $\epsilon>0$.  By Theorem 3 in \cite{LarussonSigurdsson2}, since $V_U(a)=0$, there is a map $g:\overline{\mathbb D}\to\mathbb P^n$ that extends holomorphically to a neighbourhood of $\overline{\mathbb D}$, such that $g(0)=a$, $g(\mathbb T)\subset U$, and
$$\sum_{z\in g^{-1}(H_\infty)} m_g(z)\log|z| > -\epsilon,$$
where $H_\infty$ denotes the hyperplane at infinity, which we take to be the subset of the projective space $\mathbb P^n$ with projective coordinates $[z_0:\dots:z_n]$ where $z_0=0$.  Also, $m_g(z)$ denotes the multiplicity of the intersection of $g$ with $H_\infty$ at $z$.

Note that $g^{-1}(H_\infty)$ is a finite subset of $\mathbb D$.  Let $\phi:\mathbb D\to X=\mathbb D\setminus g^{-1}(H_\infty)$ be a universal covering map with $\phi(0)=0$.  Lift $g$ to a holomorphic map $(g_0,\dots,g_n):\overline{\mathbb D}\to\mathbb C^{n+1}\setminus\{0\}$.  Then $g_0,\dots,g_n$ are bounded on $\mathbb D$, so the components $g_1/g_0,\dots,g_n/g_0$ of $g|X:X\to\mathbb C^n$ are quotients of bounded holomorphic functions, and $f=g\circ\phi:\mathbb D\to\mathbb C^n$ is a Nevanlinna disc.  Now $\phi$ has a nontangential boundary value in $\mathbb T$ at almost every point of $\mathbb T$, so $f$ has nontangential boundary values in $U$ almost everywhere on $\mathbb T$.

We need to consider the singular measure of each component $k\circ\phi$, $k=g_j/g_0$, of $f$.  We note that $k$ is a holomorphic function on a neighbourhood of $\overline{\mathbb D}$, except for finitely many poles $z$ in $g^{-1}(H_\infty)$ of (negative) order $o_k(z)\geq -m_g(z)$.  We may assume that $k$ is not identically zero (otherwise, its singular measure is zero) and that $k$ has no zeros on $\mathbb T$ (if necessary, we can slightly perturb $g$ to ensure this).  Define 
$$u(\zeta)=\sum_{z\in k^{-1}(0)} o_k(z)\log\bigg|\frac{\zeta-z}{1-\bar z\zeta}\bigg|\leq 0, \qquad \zeta\in\overline{\mathbb D}.$$
Define $v\geq 0$ by the same sum taken over the poles of $k$.  Then $u$ and $-v$ are negative subharmonic potentials on $\mathbb D$, equal to $0$ on $\mathbb T$.  Now $\log|k|-u-v$ is harmonic on $\mathbb D$, except at the zeros and poles of $k$ where it is bounded, so it extends to a bounded harmonic function $h$ on $\mathbb D$ with the same continuous boundary values as $\log|k|$.

We will now do some potential theory on the Riemann surface $X$, whose Martin boundary is $\mathbb T\cup g^{-1}(H_\infty)$.  The harmonic measure of $X$, relative to the base point $0$, is $\sigma$ on $\mathbb T$ and zero on $g^{-1}(H_\infty)$.  Thus $h$ is the Poisson integral on $X$ of the boundary function of $\log|k|$, so $h$ is the absolutely continuous part of the least harmonic majorant of $\log|k|$ on $X$.  Also, $v$ is a positive harmonic function on $X$ with zero boundary values almost everywhere, so $v$ is singular on $X$.  By definition of $h$, $\log|k|\leq h+v$, and we claim that $h+v$ is the least harmonic majorant of $\log|k|$ on $X$.  Namely, if $w\leq v$ is a harmonic function on $X$ such that $\log|k|\leq h+w$, then
$$u=\log|k|-h-v\leq w-v\leq 0$$
on $X$, so $w-v$ extends to a negative harmonic majorant of $u$ on $\mathbb D$.  Since $u$ is a potential on $\mathbb D$, we conclude that $w=v$.  Thus, $v$ is the singular part of the least harmonic majorant of $\log|k|$ on $X$.

Finally, precomposition by a covering map preserves least harmonic majorants and their decompositions, so the singular part of the least harmonic majorant of $\log|k\circ\phi|$ on $\mathbb D$ is $v\circ\phi\geq 0$.  The mass of the corresponding singular measure is
$$v(0)= \sum_{z\in k^{-1}(\infty)} o_k(z)\log|z| \leq -\sum_{z\in g^{-1}(H_\infty)} m_g(z)\log|z| < \epsilon,$$
and the proof is complete.
\end{proof}

Note that in general, the disc $f$ is no better than Nevanlinna.  To take a simple example, let $g(z)=1/z$ and $\phi(z)=\exp\dfrac{z+1}{z-1}$.  Then $g\circ\phi$ is not in the Hardy class $H^p$ for any $p>0$, nor in the class $N^+$ (see \cite{Garnett}, p.~71).  In fact, $g\circ\phi$ is the singular function $s$ whose measure $\mu$ is the unit mass at 1.

A comparison with Poletsky's characterization in \cite{Poletsky} of the polynomial hull of a pluriregular compact set is in order.  Here is a version of his result, easily derived from his fundamental theorem on the plurisubharmonicity of the Poisson envelope.  It is a slight modification of Theorem 7.4 in \cite{LarussonSigurdsson1}, with a stronger condition (ii) but an identical proof.  Recall that an analytic disc is said to be closed if it extends holomorphically to a neighbourhood of $\overline{\mathbb D}$.

\begin{theorem}[Poletsky 1993]  Let $K\subset\mathbb C^n$ be compact, $a\in\mathbb C^n$, and $\Omega$ be a pseudoconvex neighbourhood of $K$ and $a$, bounded and Runge.  Then the following are equivalent. 
\begin{enumerate}
\item[(i)] $a$ is in the polynomial hull of $K$.
\item[(ii)] For every neighbourhood $U$ of $K$ and $\epsilon>0$, there is a closed analytic disc $f$ in $\Omega$ with $f(0)=a$ and $\sigma(\mathbb T\setminus f^{-1}(U))<\epsilon$. 
\end{enumerate}
\label{poletsky}
\end{theorem}

This result is clearly stronger than ours in that it does not require $K$ to be connected.  When $K$ is connected, however, condition (ii) in Theorem \ref{characterization} implies condition (ii) in Theorem \ref{poletsky}, apart from the property that all discs lie in $\Omega$.  Namely, let $f:\mathbb D\to\mathbb C^n$ be a Nevanlinna disc with $f(0)=a$ and almost all nontangential boundary values in the open neighbourhood $U$.  (The singular measures of the components of $f$ are irrelevant here.)  Then the closed analytic discs $f_r$ with $f_r(z)=f(rz)$, $r<1$, satisfy $\sigma(\mathbb T\setminus f_r^{-1}(U))\to 0$ as $r\to 1^-$, because $\chi_U\circ f_r\to\chi_U\circ f^*$ pointwise almost everywhere on $\mathbb T$ (for this we need $U$ to be open), so
$$\sigma(\mathbb T\cap f_r^{-1}(U)) = \int_\mathbb T \chi_U\circ f_r\,d\sigma \to \int_\mathbb T \chi_U\circ f^*\,d\sigma = 1.$$
As for the converse, it is generally impossible to pass to any useful limit of analytic discs as $\epsilon\to 0$ in Theorem \ref{poletsky}, so we see no way of deriving condition (ii) in Theorem \ref{characterization} from condition (ii) in Theorem \ref{poletsky}.  Thus the two theorems differ significantly and complement each other.

\end{document}